\documentclass{article}
\usepackage{mathptmx}
\usepackage{tikz}
\usetikzlibrary{positioning,quotes}
\def\bq{\ensuremath{\mathbf{Q}}}

\def\aut(#1/#2){\mathrm{Aut}\,(#1/#2)}
\def\gal(#1/#2){\mathrm{Gal}\,(#1/#2)}
\newtheorem{lem}{Lemma}

\newtheorem{defn}{Definition}

\title{Inflated $G$-Extensions for Algebraic Number Fields}
\author{M Krithika,
VIT University, Chennai,\\
and\\
P Vanchinathan\\ 
\texttt{krithika.m2020@vitstudent.ac.in\qquad vanchinathan@gmail.com}}

\begin{document}
\maketitle

\begin{abstract}
	In 2018, Legrand and Paran proved 
	a weaker form of the Inverse Galois Problem
	for all Hilbertian fields and all finite groups
	: that is, there exist possibly non-Galois extensions over given Hilbertian base field with given finite group as the group of field automorphisms fixing the base field. For $\bq$ it was proved earlier by M. Fried.  In this paper our objective is
	how big the degree of such extension can be compared to the order of the automorphism group.
	A special case of our result shows that if the Inverse Galois problem for $\bq$ 
	has a solution for
	a finite group $G$, say of order $n$, then there exist algebraic number  fields
	of degree  $nm$, for any $m\ge3$ with the same automorphism group $G$.
\end{abstract}

\section{Introduction}
In field theory, two opposite extreme examples of extension fields are well-known. By extreme we mean in terms of the number of automorphism of an extension $K$ fixing the base field $F$ element-wise. For any Galois extension $K/F$ the order of the automorphism group is equal to the degree of the extension $K$ over $F$ which is the maximum possible. On the other hand $K=\bq(2^{1/n})$, with $n$  an odd number, is an extension of  degree $n$ over $\bq$, with the barest minimum  number of
automorphisms, namely just the trivial automorphism.  These two  extremes apart we note that
$K=\bq(2^{1/4})$  is of degree $4$ over $\bq$ with exactly two automorphisms.
The focus of our investigation being the relationship  between the degree of $K/F$ and 
and  the order of $\aut(K/F)$, 
we are naturally led to the following definition.
\begin{defn}
	For a finite separable extension $K/F$, the ratio $[K:F]/|\aut(K/F)|$ 
is called the inflation index. When a finite group  $G$ is the
	group of automorphisms of a field $K$ over $F$ we say 
	$K$ is an inflated $G$-extension of $F$. 
\end{defn}
Here are some simple properties of this inflation index:
\begin{enumerate}
	\item The inflation index is a positive integer.
		(This  follows from Lemma~\ref{NHmodH} in Section 2 below.)
\item 
If $K/F$ is a Galois extension, then the inflation index is $1$. 
		(By abuse of terminology we call this an 
		uninflated extension.)
	\item Inflation index of $K$ over $F$ is equal to the number of \textit{distinct}
		extension fields of $F$ 
inside an algebraic closure of  $ {F}$ 
		that are isomorphic to $K$ over $F$.
\item The higher the inflation index the poorer
the extension is with respect to the number of automorphisms.
\end{enumerate}
In some sense inflation index is a numerical measure of the degree of
departure from being a Galois extension. 
The extension  $\bq(2^{1/n})$ of $\bq$ mentioned above, under our new terminology, is an inflated $G$-extension with
\begin{enumerate}
\item[(i)]  inflation index $n$ for $G$ trivial,  in case  $n$ is odd.
\item[(ii)] inflation index $n/2$ for $G$ the cyclic group of order $2$,  in case $n$ is even.
\end{enumerate}

The celebrated unsolved inverse Galois problem can be stated as asking if uninflated 
$G$-extensions exist over a given field $F$. 
When we ask the weaker question of existence of inflated $G$-extensions, a result of Legrand and Paran \cite{leg} gives affirmative answer 
over all Hilbertian fields $F$ and all finite groups $G$. There has been solutions
for other cases by earlier authors to this problem. See \cite{kol}, \cite{fried},\cite{geyer} and \cite{tak}. As their aim was on realizing every finite group as 
the automorphism group, the question 
of comparing the degree of the extension to the order of the group 
was not considered.
But looking at the construction of Legrand and Paran it is not hard to see that their solution has inflation index $3k$ where $k$ is the index of $G$ as an  embedded subgroup of some $S_n$.

In our earlier work \cite{clusters} we concentrated on finding the number of roots of an irreducible polynomial that belong to the field generated by a single root, a study initiated by  
Perlis \cite{per}.
(His work  seems  to have gone largely unnoticed). We had shown for all number fields $F$
and numbers $n,r$ with $r|n$ (except the case $n=2r$), there exists an algebraic number 
$\alpha$ of relative degree $n$ over $F$, with exactly $r$ of its $F$-conjugates 
lying in  $F(\alpha)$.
This number $r$ is actually the order of the group of automorphisms of
	$F(\alpha)$ over $F$.
After becoming aware of the work of Legrand and Paran few months back, 
we were led to 
pursue the following question:

\vspace{6pt}
\noindent \framebox{\vbox{\hsize=\textwidth \noindent \it For a given finite group $G$ of order $n$, and a number field $F$,
	do there  exist inflated $G$-extensions over $F$  for any given inflation index?}}

\vspace{6pt}
Combining  their result   with  the method used in our work (\cite{clusters})  
we have now arrived at a positive answer to this  question for inflation index $\ge4$, 
and a little more.
The precise result is below:

\vspace{6pt}
\noindent\textbf{Main Theorem (Inflating an Uninflated or a Pre-inflated Extension)}:\quad
Let $F$ be an algebraic number field and $G$ be a finite group of order $n$. 
Assume there exists a field $K$ that is an inflated $G$-extension of $F$ with
inflation index $k\geq 1$. That is, $\aut(K/F)\cong G$ and $[K:F]= k\times|G|=kn$.
Then we have
\begin{enumerate}

\item[(a)]   for any given $m\geq 4$ there exist infinitely many inflated $G$-extensions with inflation index $mk$.
\item[(b)]  in case $F=\bq$,  there also exist infinitely many inflated $G$-extensions with inflation index $3k$.
\item[(c)]  in case $G$ is an abelian group  and  $F=\bq$,
	there also exist infinitely many inflated $G$-extensions with inflation index $2k$.

\end{enumerate}

This paper is organized as follows. In Section 2 we collect all the preliminaries needed to prove
the main theorem. Section 3 gives the proof of our Main Theorem spread over three subsections.
The final section provides illustrative examples. 

We want to emphasize that proof of our Main Theorem needs no previous knowledge of our paper mentioned above or other papers cited here.  Except for using a result of Hilbert from 19th century
on the existence of $S_n$ and $A_n$ Galois extensions of number fields our solution needs no prerequisites.

\vspace{6pt}
\noindent \textit{Acknowledgement:}\quad 
This research work was carried out while the 
second author was visiting the Institute of Mathematical Sciences, Chennai,  India, which provided
an excellent atmosphere to work.

%%%%%%%%%%%%%%%%%%%%%%%%%%%%%%%%

\section{Preliminaries}
As we deal with inflated extensions the following folklore lemma in field theory needs to be
singled out as it will be used repeatedly. This is Lemma 2.1 in \cite{leg} and in \cite{per}.
\begin{lem}\label{NHmodH}
Let $L/F$ be a finite separable extension with $\tilde L/F$ any Galois extension 
	of $F$ containing $L$	(not necessarily the Galois closure of $L$)  with Galois group
$G$. Then $\aut(L/F)\cong N_G(H)/H$ where $H\subset G$ is the subgroup with $L$ as its fixed field. 
Further, the  number of distinct conjugates of $L/F$ is the index of $N_G(H)$ in $G$.
\end{lem}

It is well-known that for the compositum of two linearly disjoint Galois extensions, the Galois
group is the direct product of those individual Galois groups. The following is a very useful 
generalization for automorphism groups when  the extensions under discussion  are not Galois.

\begin{lem}\label{disjcomp} 
	Let $E_1,E_2$ be any two Galois extensions over a field $F$ which are linearly disjoint  
	over $F$.
Let their  Galois groups be $G_1, G_2$ respectively with  $H_1, H_2$ their subgroups.
Let $K_1\subset E_1$ and $K_2\subset E_2$ be the fixed fields associated to $H_1, H_2$ respectively. Let $N_1, N_2$ be normalizers of $H_1, H_2$ respectively in $G_1, G_2$.
Then for the compositum $K_1K_2$  we have
\begin{itemize}
\item[(i)]  $K_1K_2$ is the fixed field of $H_1\oplus H_2$
\item [(ii)] 
$\aut(K_1K_2/F)\cong  \aut(K_1/F)\oplus \aut(K_2/F)$.
\end{itemize}
\end{lem}
Proof:\quad 
(i) is obvious.
Statement (ii) can be proved by  applying   Lemma~\ref{NHmodH} for $L=E_1E_2$, 
by noting   $\gal(E_1E_2/F)\cong G_1\oplus G_2$
and
$$N_{G_1\oplus G_2} (H_1\oplus H_2) = N_{G_1}(H_1) \oplus N_{G_2}(H_2)$$
 
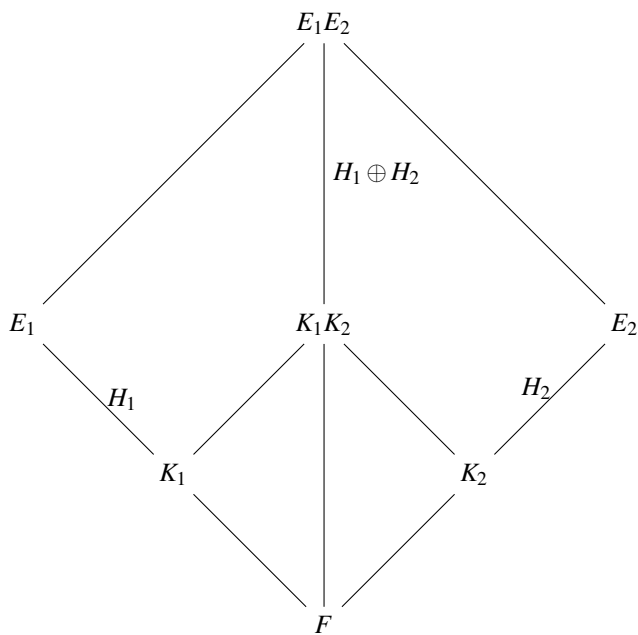
\begin{center}
\begin{tikzpicture}
    \node (Q) at (0,0) {$F$};
    \node (Q1) at (-2,2) {$K_{1}$};
    \node (Q2) at (2,2) {$K_{2}$};
    \node (Q3) at (-4,4) {$E_1$};
    \node (Q4) at (0,4) {$K_{1}K_{2}$};
    \node (Q5) at (0,8) {$E_{1}E_{2}$};
    \node (Q6) at (4,4) {$E_{2}$};

	\draw (Q)--(Q1);% node [left, pos=0.5]{} ;
	\draw (Q)--(Q2);% node [right, pos=0.5]  ;
	\draw (Q1)--(Q3)node [right, pos=0.5]{$H_1$} ;
        \draw (Q1)--(Q4);
        \draw (Q2)--(Q4);
        \draw (Q)--(Q4);% node [left, pos=0.6];
        \draw (Q3)--(Q5);
	\draw (Q4)--(Q5) node [right, pos=0.5]{$H_1 \oplus H_2$} ;
	\draw (Q2)--(Q6)node[left,pos=0.6]{$H_2$} ;
        \draw (Q6)--(Q5);
   % \draw (Q4) to[bend right=40, looseness=1,"$H_1\oplus H_2$"] (Q5); 
   %\draw (Q2) to[bend left=30, looseness=1, "$H_2$"](Q6) ; 
   %\draw (Q1) to[bend left=40, looseness=1,"$H1$"] (Q3); 
    %\draw[bend right=40]  (Q) to node anchor=west]{$H$} (Q1);
    
\end{tikzpicture}

\end{center}

Our Main Theorem stated for a general finite group is very easy for the case
the group is trivial, and possibly is well-known without our jargon. Nevertheless
considering the crucial role it plays in this work we state and prove it here:
\begin{lem}\label{trivialgroup}
	For every positive integer $m\ge3$ and any number field $F$ there exist 
	infinitely many extensions $K$ 
	of degree $m$ over $F$  with $\aut(K/F)$ trivial.
	That is, inflated extensions for trivial group exists over all algebraic  number
	fields and for all inflation  indices $\ge3$.
\end{lem}

Proof: \quad 
\textit{Clearly $m=2$ needs to be avoided as quadratic extensions for number fields
are always uninflated (that is,  Galois).}

This is readily found inside the Galois extensions $L_m/F$
with Galois group the alternating group $A_m$ (or the symmetric group  $S_m$)
which by results of Hilbert, are infinitely many.
Now look at the subgroup $A_{m-1}$ (or respectively $S_{m-1}$) embedded
as isotropy of one element for the natural action of $A_m$ (or $S_m$) on $m$ letters.
The fixed field $K_m$ for such a subgroup will do. As $A_{m-1}$ (or $S_{m-1}$) is a maximal subgroup
of index $m$ in $A_m$ (respectively in $S_m$) its normalizer is seen to be itself.
Now we can apply Lemma~\ref{NHmodH} and conclude that the only automorphism of $F_m$ over $F$ is the identity map.
(Note actually that in this case we can say a little more: there are no intermediate fields
between $F_m$ and $F$).

\section{Proof of Main Results}

\subsection{Proof of Part (a):}

First we prove under the additional hypothesis that $K/F$ is Galois. 
	So we are in the cases with $m\ge 4$, $k=1$. Let us consider the Galois extension 
	$L_m$  over $F$ guaranteed by Hilbert
	with the alternating group $A_m$ as its Galois group. 
We claim this can be chosen to be linearly disjoint with $K$ over $F$. For this purpose
it suffices to show $L_m\cap K=F$, as $L_m$ is Galois over $F$.
	When $m\ge 5$, as $A_m$ is simple  no intermediate field in $L_m$ could be Galois. 
	Consequently   the intersection
	$L_m\cap  K$  can only 
	be either $L_m$ or $F$.  For the case $m=4$, as $A_4$ admits
	a normal subgroup of index 3, and intersection could be that cubic Galois extension.
	The availability of infinitely many linearly disjoint  $A_m$-extensions of $F$ 
	ensures we can choose  $L_m$ to be  linearly disjoint with $K$.
	Now let $F_m\subset L_m$ be the fixed field for the subgroup
	$A_{m-1}\subset A_m$. This being a maximal subgroup, its normalizer is itself
	and now it follows from  Lemma~\ref{NHmodH}  that $\aut(L_{m}/F)$ is trivial.

Now consider the compositum $K_m:=KF_{m}$ and regard it  as a subfield of
the bigger compositum $KL_m$.  This $K_m$ is the fixed field for the subgroup
$A_{m-1}\oplus \{e\}$ of $A_{m}\oplus G$, in the Galois extension $KL_m/F$.
By appealing to Lemma~\ref{disjcomp} we can now see that $\aut(K_m)\cong (A_{m-1}\oplus G)/(A_{m-1}\oplus \{e\}) \cong G$, Now by linear disjointness again the degree of the extension  $K_m$ over $F$ is the product of degrees of $F_m$ and $K$ which is  $mn$ as required. Thus $K_{m}/F$ is an
inflated $G$-extension  with inflation index $m$.

This proves part (a) of our Main Theorem in the special case when $K/F$ is a Galois extension.
Now we go to the general case with inflation index $k>1$, 
So $K$ is \textit{not} Galois over $F$. Let $\tilde K$ be its Galois closure over $F$.  We apply Lemma~\ref{disjcomp} with $E_1=L_m, E_2=\tilde K$
$H_1=A_{m-1}\cong \gal(L_m/F_m),\ H_2=\gal(\tilde K/K)$ which yields
$$\aut(KF_m/F)\cong \aut(F_m/F)
\oplus \aut(K/F) = \{e\}\oplus G\cong G,$$ and the degree $KF_m=mk|G|=mkn$, completing the proof of part (a) in all cases.

\subsection{Proof of Part (b), $m=3$}
As in the proof of part (a) it suffices to prove only for the case
$K/\bq$ is Galois, i.e. $k=1$.
By hypothesis we are given a Galois  extension $K/\bq$ with Galois group $G$ for some finite
group $G$ of order $n$. We will show how to construct  extensions $K_3$  over $\bq$ of degree $3n$ with $\aut(K_3/F)\cong G$.

Let the discriminant of the Galois extension $K/\bq$ be $\Delta$.

\textbf{Case (i)  $3 \not|\Delta$}:\quad  Then choose a prime $p$ not dividing $\Delta$.
In the cubic extension $F_{3}/\bq$  given by a root of $x^3-p$, the only ramified primes are
3 and $p$ and so its Galois closure $L_3 $  and $K$ are linearly disjoint.
As $F_{3}/\bq$ has no non-trivial automorphism  we see that
for   the compositum $K_{3}=KF_{3}$ has degree $3n$ and its automorphism group is $\{e\}\oplus G$ 
as subgroup of the Galois group $S_3\oplus G$, and so isomorphic to $G$. So $K_3/F$ is an inflated $G$-extension with inflation index $3$.

\textbf{Case (ii)   $3 \mid \Delta$}:\quad  For primes $p>\Delta$  we  consider the
two Eisenstein polynomials $f(x)=x^3\pm p\Delta x+ p\Delta$ whose discriminants are
$\mp4p^3\Delta^3-27p^2\Delta^2=\mp p^2\Delta^2(4p\Delta \pm 27)$. One of these two 
 polynomials will have negative discriminant  (hence a non-square in $\bq$)  depending 
 on the sign of $\Delta$. For that choice the  Galois group will be  $S_3.$ 
Let $F_{3}/\bq$ be obtained by adjoining one root of this cubic.  We need $L_{3}$ to be linearly 
disjoint with $K$ over $\bq$. But $K \cap L_{3}$ should be a Galois extension. 
So the  possibilities for this intersection are
$\bq, L_{3}$, or its unique quadratic subfield.
To show linear disjointness it suffices to show that the above intersection does not contain
the unique quadratic subfield of $L_3$.

That  quadratic extension inside $L_{3}/\bq$ is obtained by adjoining the
square root of $4p\Delta-27$. We need to choose a prime $p>\Delta$ such that $4p\Delta-27\neq 3^k$ for some $k\geq 0$ then there exist a new prime which divides the discriminant and not dividing $\Delta$ will not ramify in $K$ so they are linearly disjoint.  It suffices to show for large $p$ $4p\Delta-27$ is not a power of 3.

We claim  that such primes  are at most finitely many.  So assume  $4p\Delta +27 =3^k$ for some $k$, equivalently $4p\Delta =3^k-27$ ($k=0,1,2,3$ is not possible). We have $4p\Delta=27(3^{(k-3)}-1)$ since $\gcd(27,4p)=1$ 
implies $\Delta $ to be a multiple of $27$, but not a multiple of
higher power of $3$. We denote by $\Delta':=\Delta/27$, an integer coprime to 3.
We have  $4p\Delta'=3^{k-3}-1$ and so $4p\Delta'\equiv1\pmod{3}$ so $p\Delta' \equiv1\pmod{3}$. Since $\Delta'$ is fixed this cannot be true for all primes, proving our claim and thereby giving
us linear disjointness as desired.

Hence the compositum $K_3=KF_{3}$ has degree $3n$  and its automorphism group is $\{e\}\oplus G$ as subgroup of the Galois group $S_3\oplus G$ and is isomorphic to $G$.

\subsection{Proof of Part (c), $m=2$}
Let $C_n$ denote the cyclic group of order $n$ written multiplicatively with $g$ a generator. Let $D_n$ denote the dihedral group of order $2n$ 
whose elements may be listed conventionally as below:
$$D_n=\{ 1,x,x^2,\cdots,x^{n-1},y,yx,\cdots,yx^{n-1}\}$$

\begin{lem}
There exist Galois extension of $\bq$ with Galois group $C_n\times D_n$.
\end{lem}
\textbf{Proof:}
We assume the well-known result that Galois extension over $\bq$ with $D_n$ as Galois group exists.
Let $K_1$ be such an extension, with discriminant $\Delta$.

Any $p^{\rm th}$-cyclotomic extension with $p$ not dividing $\Delta$ will be linearly disjoint with $K_1$.
Choose $p$ such that $p\equiv 1 \pmod {n}$ and $p>\Delta$ this is possible (by Dirichlet's Theorem on arithmetic progression). This $p^{th}$ cyclotomic extension will contain a cyclic extension $K_2$ of degree $n$ over $\bq$ again linearly disjoint over $\bq$. Now the compositum $K=K_{1}K_2$ over $\bq$ has Galois group $D_n\times C_n$.

Before proving part (c) of our Main theorem for general abelian group
we do it for the special case of cyclic groups as below:
\begin{lem}
	With $K$ as above a Galois extension of $\bq$ with $D_n\times C_n$ as Galois group for $H\subset D_n\times C_n$ defined as $H=\{ (x^i,g^i) \mid i=0,1,2,\cdots,(n-1)\}$ the fixed field $L$ of $H$, $L\subset K$ is an inflated $C_n$-extension with inflation index $2$.
\end{lem}
\textbf{Proof:}

Given $H=\{(x^i,g^i) \mid i=0,1,2,\cdots,(n-1)\}$ we can see that 
$N(H)=\{(x^i,g^j)\mid i,j=0,1,2,\cdots,(n-1)\}$ so $N(H)\cong C_n \oplus C_n\subset D_n \oplus C_n$
Thus $N(H)/H \cong C_n$ and the index of $H$ in $D_n\oplus C_n$ is $2n^2/n=2n$. So the degree of the fixed field for $H$ is twice the order of its automorphism group thereby 
	proving the Lemma.

Let $G$ be any finite abelian group we now 
need to show inflated $G$-extension with inflated index $2$ exists.
 This will be done by some kind of induction on the number of cyclic factors in the decomposition
of the abelian group. The Lemma proved just now  starts the induction.
Using the structure theorem  $G$ can be expressed as  $G\cong C_r\oplus A$, 
for some cyclic group $C_r$ and some  abelian group $A$. 

Let $K/\bq$ be an extension of degree $2r$ with automorphism group cyclic of order $r$ just constructed above, with $\widetilde{K}$ its Galois closure. Now for the abelian group $A$ again using Dirichlet theorem, we can find Galois extension $L/\bq$ with $\gal(L/\bq)\cong A$ and $L$ linearly disjoint with $\widetilde{K}$. Then for the compositum,
$$\aut(KL/\bq)\cong \aut(K/\bq) \oplus \aut(L/\bq)\cong C_r \oplus A$$ 
And the degree is computed as twice the order of $G$.
$$[KL:\bq]=[K:\bq][L:\bq]=2r\times |A|=2|C_r\oplus A|=2|G|$$
Thus we have proved all the parts of our Main Theorem.

\section{Examples}

\noindent
\textbf{Inflated $C_3$-extensions over $\bq$ with inflation index 4}

	\vspace{6pt}
Note that the discriminant of the polynomial $f(x)=x^3-3x+1$ is $81$, a perfect
square, whose root yield a cyclic cubic  extension $L/\bq$.
  The only prime ramifying  here is $3$. We construct 
	an $A_4$-extension  $K/\bq$ from another  cubic extension 
	embedded in $\bq(\zeta_p+\zeta_p^{-1})$
	for $p\equiv1\pmod 3$. By Dirichlet's theorem we have infinitely many choices 
	for $p$ allowing us to pick one that does not divide the discriminant 
	of  $L$. The above process guarantees the linear disjointness. Specifically
	take $p=7$. Suitable quadratic extension of the real subfield of
	the 7th cyclotomic field will be non-Galois, with their Galois closure
	yielding $A_4$-extensions over $\bq$. Computing with the SAGE software package gives us
	$\bq(\sqrt{ \zeta_7 +\zeta_7^6},\sqrt{\zeta_7^2+\zeta_7^5},\sqrt{\zeta_7^3+\zeta_7^4})$ has degree $12$ over $\mathbf{Q}$  with $A_4$ as its Galois group. The degree $4$ extension $E$ inside this degree $12$ extension is generated by a root of
	the quartic irreducible polynomial $x^4+8x^2+64x+144$ whose discriminant is $2^67^2$ therefore the fields $K$ and $E$ are linearly disjoint over $\bq$.
The compositum $EK$ is generated by a root of  $x^{12}+12x^{10}+188x^9+534x^8+3108x^7+23738x^6+73860x^5+266037x^4+1053896x^3+2414142x^2+2634348x+2870297$ this will be a $C_3$-extension of inflation index $4$.

	\vspace{9pt}
\noindent \textbf{Inflated $S_3$-extensions over $\bq$  with inflation index $3$}

	\vspace{6pt}
Again our computations were carried out with SAGE.
For $f(x)=x^3-x+1$ the discriminant  is $-23$ which is a non-square thus it is not Galois take the splitting field $K$ of $f(x)$ which is generated by $x^6 - 6x^4 + 9x^2 + 23$ the Galois group of this polynomial is $S_3$. Thus $K/\mathbf{Q}$ is a degree $6$ Galois extension with $S_3$ as its Galois group. For a prime $p\neq 23$ the polynomial  $x^3-p$ has  discriminant is $p^23^3$ the only primes ramifying here are $p$ and $3$ so the splitting field $\mathbf{Q}(\zeta_3,\sqrt[3]{p})$ which is linearly disjoint with $K$ so the compositum of $\mathbf{Q}(\sqrt[3]{p})$ and $K$ is a degree $18$ extension over $\bq$. When $p=5$ we have the polynomial $x^{18}-18x^{16}-30x^{15}+135x^{14}+180x^{13}-96x^{12}-270x^{11}+1287x^{10}+15830x^9+4293x^8-124650x^7+275814x^6+133380x^5-284958x^4+2970x^3+446958x^2+358830x+188217$ generates an inflated  $S_3$-extension with inflation index $3$.

	\vspace{9pt}

\noindent\textbf{Tower of $C_5$-extensions with inflation index $5^n$.}

	\vspace{6pt}
Just as a curiosity we mention this example.
Let us start with  any $C_5$-extension $K_1/\bq$. There exist infinitely many
extensions $L_{1}/\bq$ with
trivial automorphism group and with inflation index $5$, as a consequence of
Lemma \ref{trivialgroup}. Choose one which is linearly disjoint with $K_{1}$ take the compositum with $L_{1}$ call it as $K_{2}$ this is an inflated
$C_5$-extension of $\bq$ with inflation index $5$. 
Now use the fact that Galois extensions over number fields with $S_n$ as Galois groups
can be chosen to be linearly disjoint with any desired extension by choosing the
ramifying primes to be coprime to the discriminant of that given field.

In the next step take this $K_2$ as base field and construct an inflated  $C_5$-extension over
$K_2$ with inflation index 5. And repeat this process obtaining $K_{n+1}/K_n$ an inflated
$C_5$-extension of inflation index 5.

So in this infinite tower 
$$\bq\subset K_1\subset K2 \subset \ldots \subset K_n\subset \ldots$$
regarding all of them as extensions of $\bq$ we get inflated $C_5$-extensions $K_{n+1}$ 
of inflation index  $5^n$.

\end{document}